\documentclass[reqno,12pt]{amsart}
\usepackage{amscd,amssymb}
 
\theoremstyle{plain}
\newtheorem{Thm}[subsection]{Theorem}
\newtheorem{Cor}[subsection]{Corollary}
\newtheorem{Lem}[subsection]{Lemma}
\newtheorem{Prop}[subsection]{Proposition}

\theoremstyle{definition}
\newtheorem{Def}[subsection]{Definition}

\theoremstyle{remark}

\newtheorem{Rem}[subsection]{Remark}

\numberwithin{equation}{section}
\renewcommand{\rm}{\normalshape}

\newcommand{\on}{\operatorname}

\newif\ifShowLabels
\ShowLabelstrue
\newdimen\theight
\def\TeXref#1{%
	\leavevmode\vadjust{\setbox0=\hbox{{\tt
		\quad\quad  {\small \rm #1}}}%
	\theight=\ht0
	\advance\theight by \lineskip
	\kern -\theight \vbox to
	\theight{\rightline{\rlap{\box0}}%
	\vss}%
	}}%

\ShowLabelsfalse

\renewcommand{\sec}[2]{\section{#2}\label{S:#1}%
	\ifShowLabels \TeXref{{S:#1}} \fi}
\newcommand{\ssec}[2]{\subsection{#2}\label{SS:#1}%
	\ifShowLabels \TeXref{{SS:#1}} \fi}

\newcommand{\reft}[1]{Theorem ~\ref{T:#1}}
\newcommand{\refl}[1]{Lemma ~\ref{L:#1}}

\newcommand{\refe}[1]{\eqref{E:#1}}

\newenvironment{thm}[1]%
	{ \begin{Thm} \label{T:#1}  \ifShowLabels \TeXref{T:#1} \fi }%
	{ \end{Thm} }

\renewcommand{\th}[1]{\begin{thm}{#1} \sl }
\renewcommand{\eth}{\end{thm} }

\newenvironment{lemma}[1]%
	{ \begin{Lem} \label{L:#1}  \ifShowLabels \TeXref{L:#1} \fi }%
	{ \end{Lem} }
\newcommand{\lem}[1]{\begin{lemma}{#1} \sl}
\newcommand{\elem}{\end{lemma}}

\newenvironment{propos}[1]%
	{ \begin{Prop} \label{P:#1}  \ifShowLabels \TeXref{P:#1} \fi }%
	{ \end{Prop} }
\newcommand{\prop}[1]{\begin{propos}{#1}\sl }
\newcommand{\eprop}{\end{propos}}

\newenvironment{corol}[1]%
	{ \begin{Cor} \label{C:#1}  \ifShowLabels \TeXref{C:#1} \fi }%
	{ \end{Cor} }
\newcommand{\cor}[1]{\begin{corol}{#1} \sl }
\newcommand{\ecor}{\end{corol}}

\newenvironment{defeni}[1]%
	{ \begin{Def} \label{D:#1}  \ifShowLabels \TeXref{D:#1} \fi }%
	{ \end{Def} }
\newcommand{\defe}[1]{\begin{defeni}{#1} \sl }
\newcommand{\edefe}{\end{defeni}}

\newenvironment{remark}[1]%
	{ \begin{Rem} \label{R:#1}  \ifShowLabels \TeXref{R:#1} \fi }%
	{ \end{Rem} }
\newcommand{\rem}[1]{\begin{remark}{#1}}
\newcommand{\erem}{\end{remark}}

\newcommand{\eq}[1]%
	{ \ifShowLabels \TeXref{E:#1} \fi 
	   \begin{equation} \label{E:#1} }
\newcommand{\eeq}{ \end{equation} }

\newcommand{\prf}{ \begin{proof} }
\newcommand{\epr}{ \end{proof} }

\newcommand\sig{\sigma}

\newcommand\calK{{\mathcal{K}}}
\newcommand\calL{{\mathcal{L}}}

\newcommand\calN{{\mathcal{N}}}
\newcommand\calO{{\mathcal{O}}}
\newcommand\calP{{\mathcal{P}}}
\newcommand\calQ{{\mathcal{Q}}}

\newcommand\calS{{\mathcal{S}}}

\newcommand\calV{{\mathcal{V}}}

\newcommand\bfd{{\mathbf d}}		
		
		\newcommand\bfF{{\mathbf F}}

		\newcommand\bfP{{\mathbf P}}

\newcommand\CC{\mathbb{C}}

\newcommand\NN{\mathbb{N}}

	\newcommand\grb{{\mathfrak{b}}}

	\newcommand\grg{{\mathfrak{g}}}

\newcommand\sdp{\times \hskip -0.3em {\raise 0.3ex
\hbox{$\scriptscriptstyle |$}}} 

\newcommand\codim{\operatorname{codim}}

\newcommand\End{\operatorname{End\,}}

\newcommand\gl{{\bf gl}}

\newcommand\Perv{\operatorname{\Perv}}

\newcommand\x{\times}
\newcommand\ten{\otimes}

\def\pnd{\calP_{n,d}}

\newcommand\pind{\pi^{n,d}}
\newcommand\qnd{\calQ_{n,d}}

\newcommand\xnd{X_{n,d}}

\newcommand\IC{\text{IC}}
\newcommand\nnd{\calN_{n,d}}
\newcommand\tn{{\widetilde \nnd}}
\newcommand\Spr{\text{\bf Spr}}
\newcommand\tg{\widetilde \grg}

\renewcommand\IC{\text{IC}}
\newcommand\tgnd{\widetilde \grg_{n,d}}
\newcommand\gd{\grg_d}
\newcommand\bgnd{\widetilde \grg_{n,\bfd}}
\begin{document}

\title[On Ginzburg's lagrangian construction]{On
Ginzburg's lagrangian construction of representations of
$GL(n)$}

\author[A.~Braverman and D.~Gaitsgory]
{Alexander Braverman and Dennis Gaitsgory}

\thanks{The work of both authors was partially supported by the
National Science Foundation and by the Ellentuck Fund}
\address{2-175, Department of Mathematics\\
	 Massachusetts Institute of Technology
	 77 Massachusetts Ave., Cambridge MA, USA and Institute for
Advanced Study, Princeton, USA}
\address{Department of Mathematics
	Harvard University
	1 Oxford st. Cambridge MA, USA and Institute for Advanced
Study, Princeton, USA}

\email{braval@math.mit.edu, gaitsgde@math.harvard.edu}
\maketitle

\begin{abstract} 
In \cite{Gi} V.~Ginzburg observed that one can realize
irreducible representations of the group $GL(n,\CC)$ in the
cohomology of certain Springer's fibers for the group
$GL_d$ (for all $d\in \NN$). However, Ginzburg's construction of
the action of
$GL(n)$ on this cohomology was a bit artificial (he defined
the action of Chevalley generators of the Lie algebra $gl_n$ 
on the corresponding cohomology by certain explicit correspondences,
following the work of A.~Beilinson, G.~Lusztig and R.~MacPherson
(\cite{BLM}), who gave a similar construction of the quantum
group $U_q(gl_n)$). 

In this note we give a very simple geometric definition of the
action of the whole group $GL(n,\CC)$ on the above
cohomology and simplify considerably the results of \cite{Gi}.
\end{abstract}

\sec{}{Introduction}

\ssec{}{The combinatorics}
Throughout this note we will fix natural numbers $n$ and $d$. $E$ will be a fixed $n$-dimensional (complex) vector space $E$ and
$G$ will denote the group $GL(E)$. Let, in addition, $\grg_d$ denote the Lie algebra $\gl(d,\CC)$.

We set
\begin{align*}
\pnd=\{(n_1,...,n_k)|\ \sum\limits_{i=1}^k n_i=d;\quad n_i\leq n\ 
\text{for all $i=1,...,k$\ } \\
\text{and $n_1\geq n_2\geq ...\geq n_k\geq 1$}\}
\end{align*}

With any $\bfP\in\pnd$ one usually associates three objects: an irreducible
representation $\rho(\bfP)$ of the symmetric group $S_d$,
an irreducible representation $V(\bfP)$ of $G$ and a nilpotent orbit 
$\calO_{\bfP}$ in $\grg_d$.

First of all, $\calO_{\bfP}$ by definition consists of all nilpotent
matrices in $\grg_d$, which have $k$ Jordan blocks of sizes
$n_1,...,n_k$. 

The representation $\rho(\bfP)$ is characterized uniquely by the following property:

Let $S_{\bfP}=S_{n_1}\x ...\x S_{n_k}$ be the corresponding
subgroup of $S_d$. If the partition $\bfP$ is a (strict) refinement of a partition $\bfP'$, we shall write 
$\bfP<\bfP'$. We have:
$$\on{Hom}_{S_d}(\rho(\bfP), \on{Ind}_{S_{\bfP}}^{S_d} {\bf 1})\neq 0, \,\,
\on{Hom}_{S_d}(\rho(\bfP), \on{Ind}_{S_{\bfP'}}^{S_d} {\bf 1})=0, \text{ if } \bfP<\bfP'.$$

Finally we
set $V(\bfP)=\on{Hom}_{S_d}(\rho(\bfP),E^{\ten d})$.  

\ssec{variety}{The variety $\xnd$}Fix $n$ and $d$ as above. Let
$X_{n,d}$ be the complex algebraic variety whose 
$\CC$-points consist of all partial flags (filtrations) 
$\calV$ of the form
\eq{}
0=V_0\subset V_1\subset ...\subset V_{n-1}\subset V_n=\CC^d
\end{equation}

The variety $X_{n,d}$ splits into connected components numbered by the set $\qnd$ of $n$-tuples of non-negative integers 
$\bfd=(d_1,...,d_n)$ with $d_1+...+d_n=d$:

For $\bfd\in\qnd$ as above, the corresponding connected component (denoted as $\xnd^{\bfd}$) consists of flags $\calV$, for which 
$\dim(V_i/V_{i-1})=d_i$. We set $|\bfd|=\dim \xnd^{\bfd}$.

\medskip

Let now 
\eq{}
 \nnd=\{a\in \grg_d|\ a^n=0\}. 
\end{equation}
$\nnd$ is a closed subvariety of the nilpotent cone of $\grg_d$.

Let $\tn=T^*X_{n,d}$ denote the cotangent bundle to $X_{n,d}$.
It may be identified with the set of all pairs 
\eq{}
\{\calV\in X_{n,d}, a\in \grg_d|\ \text{such that}\ 
a(V_i)\subset V_{i-1}\}.
\end{equation}

For $\bfd\in\qnd$ we shall denote by $\tn^{\bfd}$ the corresponding connected component of $\tn$.

\medskip

Let $\pi^{n,d}:\tn\to \grg_d$ denote the map, given by
$\pind((\calV,a))=a$.
Then $\pind$ is actually a map from $\tn$ to $\nnd$
and it is a semi-small (cf. \cite{BM}) resolution of 
singularities of $\nnd$.

\ssec{}{The main result}

\bigskip

Recall that we have fixed an $n$-dimensional complex vector space
$E$ and we have set $G=GL(E)$.

Fix now a decomposition 
\eq{}
E=E_1\oplus ...\oplus E_n\quad \text{such that}\ \dim E_i=1
\end{equation}

For $\bfd\in\qnd$ we shall denote by $E^\bfd$ the $1$-dimensional vector space $E_1^{\ten d_1}\ten...\ten E_n^{\ten d_n}$,
where $\bfd=(d_1,...,d_n)$.

Let $T\subset G$ be the maximal torus corresponding to the decomposition $E=E_1\oplus ...\oplus E_n$.
An element $\bfd\in \qnd$ defines a character of $T$: this is the character by which 
$T$ acts on $E^{\bfd}$.

\medskip

Let now $\calL_E^{n,d}$ denote the perverse sheaf on $\tn$ described as follows: 

For $\bfd\in\qnd$, we set the restriction of $\calL_E^{n,d}$ to $\tn^{\bfd}$ to be the $1$-dimensional constant sheaf
$E^{\bfd}[2|\bfd|]$. We endow $\calL_E^{n,d}$ with a $T$-equivariant structure by letting $T$ act
on $E^{\bfd}[2|\bfd|]$ via the corresponding character.

\medskip

The following result is essentially due to V.~Ginzburg (cf. \cite{Gi})
in a slightly less invariant form.
The main purpose of this paper is to provide a simple proof
of this result.

\th{main}
There is a canonical $T$-equivariant isomorphism
\eq{main}
\pind_*(\calL^{n,d}_E)\simeq \bigoplus\limits_{\bfP\in \pnd}\IC_{\bfP}
\otimes V(\bfP),
\end{equation}
where $\IC_{\bfP}$ denotes the intersection
cohomology sheaf on the closure of $\calO_{\bfP}$ in $\nnd$.
\eth

\noindent
{\bf Remark.}\quad
Note that the right hand side of \refe{main} is independent of 
the decomposition $E=E_1\oplus ...\oplus E_n$, though the sheaf
$\calL^{n,d}_E$ does depend on it.

\medskip

For any $\bfd\in \qnd$ and $\bfP\in\pnd$ let us denote by $V(\bfP)_{\bfd}$ the
corresponding weight space of $V(\bfP)$. It is easy to
see that $\qnd$ exhaust all the weights that can appear in $V(\bfP)$, i.e.
\eq{}
V(\bfP)=\bigoplus\limits_{\bfd\in \qnd}V(\bfP)_{\bfd}.
\end{equation}

By decomposing both sides of \refe{main} with respect to the characters of $T$, we obtain 
 
\cor
There is a canonical $T$-equivariant isomorphism
\eq{mainmu}
\pind_*(\calL^{n,d}_E|_{\tn^{\bfd}})
\simeq \bigoplus\limits_{\bfP\in \pnd}\IC_{\bfP}
\otimes V(\bfP)_{\bfd}
\end{equation}
\ecor

By applying \reft{main} to $E=\CC^n$ with the standard decomposition
into one-dimensional subspaces 
and the definition of intersection cohomology 
sheaf we obtain the following result:

Consider
\eq{}
H^{\text{top}}((\pind)^{-1}(a),\CC):=\bigoplus\limits_{\bfd\in \qnd}
H^{\frac{1}{2}\codim_{\nnd}(\calO_{\bfP})}((\pind)^{-1}(a)\cap \xnd^\bfd,\CC)
\end{equation}

\cor{ginzburg}(cf. \cite{Gi})Let $\bfP\in\pnd$ and let $a\in\calO_{\bfP}$.
The vector space $H^{\text{top}}((\pind)^{-1}(a),\CC)$ can be naturally
identified with the space of the representation $V(\bfP)$ of $GL(n,\CC)$. Under this identification,
the direct summand
$$H^{\frac{1}{2}\codim_{\nnd}(\calO_{\bfP})}((\pind)^{-1}(a)\cap
\xnd^\bfd,\CC)\subset H^{\text{top}}((\pind)^{-1}(a),\CC)$$
goes over to the weight space $V(\bfP)_{\bfd}\subset V(\bfP)$.
\ecor

\noindent
{\bf Remark.}\quad 
The fact that $\pind$ is semi-small means that for $a\in \calO_{\bfP}$
the
variety $(\pind)^{-1}(a)$ has 
no irreducible components of dimension larger than 
$\frac{1}{2}\codim_{\nnd}(\calO_{\bfP})$.
Hence the vector space 
$$H^{\frac{1}{2}\codim_{\nnd}(\calO_{\bfP})}((\pind)^{-1}(a)\cap \xnd^\mu,\CC)$$
is endowed with a natural basis, consisting of the fundamental classes
of irreducible components of dimension 
$\frac{1}{2}\codim_{\nnd}(\calO_{\bfP})$.
\sec{}{Proof of \reft{main}}
\ssec{}{Digression on the Springer correspondence}

Let $\grg$ be any reductive Lie
algebra and let $W$ denote the corresponding Weyl group. Let us choose
an invariant non-degenerate bilinear form on $\grg$, by means of which 
we shall identify $\grg$ with its dual space $\grg^*$. 

Recall the construction of the Springer sheaf on $\grg$. Let $\tg$ denote the variety of all pairs
$(\grb,x)$, where $\grb$ is a Borel subgroup of $\grg$ and
$x\in \grb$. Let $p:~\tg\to \grg$ denote the natural morphism. 

We set $\Spr=p_*\CC_{\tg}[\dim \grg]$. 

It is well-known 
(cf. \cite{BM}) that $\Spr$ is a semi-simple perverse sheaf on $\grg$. Moreover,
one has a canonical isomorphism 
\eq{springer}
\End(\Spr)\simeq \CC[W]
\end{equation}
where $\CC[W]$ denotes the group algebra of $W$. 

\medskip

If now $\rho$ is any (not necessarily irreducible)
finite-dimensional representation of $W$, we attach to it a perverse sheaf $\calS_{\rho}$ on $\grg$ by setting
\eq{sro}
\calS_{\rho}=(\Spr\ten \rho)^W
\end{equation}

It follows easily from \refe{springer} that $\calS_{\rho}$ is semi-simple
and is irreducible if and only if $\rho$ is irreducible.

\medskip

Let $\bfF$ denote the Fourier-Deligne
transform functor, that maps the category
of $\CC^*$-equivariant perverse sheaves on $\grg$ to 
itself (this functor is
well-defined, since we have chosen an identification $\grg\simeq
\grg^*$).

We shall need the following version of the Springer correspondence for the Lie algebra $\grg=\grg_d$ (cf. \cite{BM}, \cite{CG}):

\th{springer} 

Let $\rho=\rho(\bfP)$ for $\bfP\in\pnd$. Then $\bfF(\calS_\rho)\simeq \IC_\bfP$,
where $\IC_\bfP$ is viewed as a perverse sheaf on $\grg_d$ via the closed embedding $\nnd\hookrightarrow\grg_d$.

\eth

\ssec{}{A reformulation of \reft{main}}

Consider $E^{\ten d}$ as a representation of the group $S_d\times G$. Since our constructions are functorial,
the sheaf $\calS_{E^{\ten d}}$ carries a $G$-action and hence a $T$-action.

\medskip

Using \reft{springer}, we can reformulate \reft{main} as follows:

\th{reformulation}One has a canonical $T$-equivariant isomorphism
\eq{reform}
\pi^{n,d}_*(\calL_E^{n,d})\simeq \bfF(\calS_{E^{\ten d}})
\end{equation}
(note that $T$ acts naturally on both sides of \refe{reform}).
\eth

The rest of this section is occupied by the proof of \reft{reformulation}.

\prf Since $\bfF$ is an involution on the category of $\CC^*$-equivariant
sheaves on $\grg$, we must construct an isomorphism
between $\bfF(\pi^{n,d}_*(\calL_E^{n,d}))$ and $\calS_{E^{\ten d}}$.

First, we shall give an explicit description of the sheaf $\bfF(\pi^{n,d}_*(\calL_E^{n,d}))$. 

\medskip

Let $\tgnd$ denote the variety of all pairs $(\calV,x)$ where
$\calV\in \xnd$ and $x\in \grg_d$ such that 
$x(V_i)\subset V_i$. We have a natural projections
$p:~\tgnd\to \gd$ and $q:~\tgnd\to\xnd$ given by
$p((\calV,x))=x$ and $q((\calV,x))=\calV$, respectively.

The variety $\tgnd$ is a vector bundle over $\xnd$ and we shall denote by $\tgnd^{\bfd}$ the preimage of the
corresponding connected component of $\xnd$.

It is known that $p$ is a small map in the sense of Goresky-MacPherson  and that
its restriction to any connected component of $\tgnd$ is
surjective. In particular, any connected component of
$\tgnd$ has dimension $d^2=\dim\grg_d$.

\medskip

We introduce a perverse sheaf $\calK_E$ over $\tgnd$ as follows: 

For $\bfd\in\qnd$ its direct summand (denoted by $\calK_E^{\bfd}$) that lives over the connected component $\tgnd^{\bfd}$ 
is set to be $E^{\bfd}[d^2]$.

\lem{four-iso}
One has a canonical $T$-equivariant isomorphism 
$$\bfF(\pi^{n,d}_*(\calL_E^{n,d}))\simeq p_*\calK_E.$$
\elem

\prf Both $\tn$ and $\tgnd$ are vector subbundles of the trivial
bundle $\xnd\x \grg_d$ over $\xnd$. Since we have chosen an 
identification $\grg_d\simeq \grg_d^*$, this trivial bundle is
canonically self-dual and 

\eq{perp}
\tn=\tgnd^{\perp}
\end{equation} 
where $\tgnd^{\perp}$ denotes the orthogonal
complement (in the trivial bundle $\xnd\x \grg_d$) to the vector 
subbundle $\tgnd$. 

Let ${\widetilde \bfF}$ denote the Fourier-Deligne transform functor that maps the category of $G_m$-equivariant 
perverse sheaves on $\xnd\x \grg_d$ to itself. It follows from the standard properties of the functor ${\widetilde \bfF}$
that one has a canonical $T$-equivariant isomorphism
$${\widetilde \bfF}(\calL_E^{n,d})\simeq \calK_E^{n,d}.$$ 

This observation, together with the fact that Fourier transform commutes
direct images, implies \refl{four-iso}.
\epr

\ssec{}{}

It remains to show that $p_*\calK_E$ is isomorphic
to $\calS_{E^{\ten d}}$. But the map $p$ is small in the sense
of Goresky-MacPherson (cf. \cite{BM}) and, therefore, $p_*\calK_E$
is equal to the Goresky-MacPherson extension of its restriction
on the set $\gd^{rs}$ of regular semisimple elements
in $\gd$. Since $\calS_{E^{\ten d}}$ is also equal to the
Goresky-MacPherson extension of its restriction on $\gd^{rs}$, 
it is enough to construct a $T$-equivariant isomorphism of the restrictions
of  $\calS_{E^{\ten d}}$ and $p_*K_E$ to $\gd^{rs}$.

Let $\CC^{(d)}$ denote the $d$-th symmetric power of $\CC$, i.e. the
quotient of $\CC^d$ by the natural action of the symmetric group $S_d$. 
Let $^{0}{\CC^d}$ denote the complement to all the diagonals
in $\CC^d$, i.e. 
\eq{}
^0\CC^d=\{ (z_1,...,z_d)\in \CC^d|\ z_i\neq z_j\ \text{for }i\neq j\}
\end{equation}
Set $^0\CC^{(d)}$ to be the image of $^0\CC^d$ in 
$\CC^{(d)}$ and we denote by $\sig_d$ the natural map
$^0\CC^d\to {^0\CC^{(d)}}$. 

For $\bfd\in\qnd$ let  
$^0\CC^{(\bfd)}$ denote the quotient of
$^0\CC^d$ by $S_{\bfd}$. Let
$\sig_{\bfd}:{^0\CC^{(\bfd)}}\to
{^0\CC^{(d)}}$ be the natural map. 

One has a Cartesian square
\eq{}
\begin{CD}
\bgnd^{rs} @>{e_{\bfd}}>> ^0\CC^{(\bfd)}\\
@VpVV               @V{s_{\bfd}}VV\\
\gd^{rs} @>{e_d}>>        ^0\CC^{(d)}
\end{CD}
\end{equation}
Here $e_d$ is the map which associates to any $x\in \gd^{rs}$ the
collection of its eigenvalues and $e_{\bfd}$ is the map which
associates to every pair $(\calV,x)\in \bgnd$ the eigenvalues of
$x$ in the spaces $V_1,V_2/V_1,...,V_n/V_{n-1}$ respectively. 

\medskip

Let $\on{S}_E$ denote the $d$-th symmetric power of the constant local system with fiber $E$ on $\CC$. In other words,
$\on{S}_E$ is a locally constant sheaf on ${^0\CC^{(d)}}$ defined as
$$\on{S}_E:=(E^{\ten d}\ten (\sig_d)_*\CC)^{S_d}.$$

We have a $T$-equivariant identification:
$$e_d^* \on{S}_E[d^2]\simeq \calS_{E^{\ten d}}.$$

\smallskip

Analogously, for $\bfd\in\qnd$, let $\on{K}^{\bfd}_E$ be the constant sheaf over $^0\CC^{(\bfd)}$ with fiber $E^\bfd$. We have:
$$e_{\bfd}^* \on{K}^{\bfd}_E[d^2]\simeq \calK_E$$
(in a $T$-equivariant way).

\medskip

The assertion of \reft{reformulation} now follows from the following lemma.

\lem{} There exists canonical $T$-equivariant isomorphism 
\eq{iso-bred}
\bigoplus\limits_{\bfd\in \qnd}(\sig_{\bfd})_*\on{K}^{\bfd}_E\simeq \on{S}_E
\end{equation}
\elem

\prf

The lemma follows easily from the following linear algebra observation:
there exists canonical $T\times S_d$-equivariant isomorphism
\eq{}
\bigoplus\limits_{\bfd\in \qnd}E^{\bfd}\simeq E^{\ten d}
\end{equation}
The proof is clear.
\epr
\epr

\ssec{}{Acknowledgments}We thank V.~Ginzburg for explaining to us the
contents of \cite{Gi}.

\end{document}